\documentstyle[12pt]{article}
\newlength{\abstractwidth}
\renewcommand{\thefootnote}{\fnsymbol{footnote}}
\renewcommand{\thanks}[1]{\footnote{#1}} 
\newcommand{\starttext}{
\setcounter{footnote}{0}
\renewcommand{\thefootnote}{\arabic{footnote}}}

\newcommand{\be}{\begin{equation}}
\newcommand{\bea}{\begin{eqnarray}}
\newcommand{\eea}{\end{eqnarray}}
\newcommand{\ee}{\end{equation}}

\def\ba{\begin{eqnarray}}
\def\ea{\end{eqnarray}}

\def\v{\vskip .1in}

\def\al{\alpha}
\def\b{\beta}

\def\g{\gamma}
\def\l{\lambda}

\def\o{\omega}

\def\ti\tilde

\def\u{\underline}

\def\pl{\partial}

\def\ddb{\partial\bar\partial}

\def\det{{\rm det}}

\def\ti{\tilde}

\def\o{\omega}

\def\[{{\bf [}}
\def\]{{\bf ]}}

\def\pl{\partial}

\v

\begin{document} \starttext \baselineskip=15pt
\setcounter{footnote}{0} \newtheorem{theorem}{Theorem}
\newtheorem{lemma}{Lemma} \newtheorem{corollary}{Corollary}
\newtheorem{definition}{Definition}
\begin{center} {\Large \bf MULTIPLIER IDEAL SHEAVES AND THE K\"AHLER-RICCI FLOW}
\footnote{Research supported in part by National Science Foundation
grants DMS-02-45371, DMS-06-04657
 and DMS-05-14003}
\\
\bigskip
{\large D.H. Phong$^{*}$, Natasa Sesum$^{*}$, and Jacob Sturm$^\dagger$}
\bigskip
\end{center}
\v
\begin{abstract}

{\small Multiplier ideal sheaves are constructed as obstructions to
the convergence of the K\"ahler-Ricci flow on Fano manifolds,
following earlier constructions of Kohn, Siu, and Nadel, and using the
recent estimates of Kolodziej and Perelman.}

\end{abstract}

\bigskip
\baselineskip=15pt \setcounter{equation}{0} \setcounter{footnote}{0}

\section{Introduction}
\setcounter{equation}{0}

The global obstruction to the existence of a Hermitian-Einstein
metric on a holomorphic vector bundle is well-known to be encoded in
a destabilizing sheaf, thanks to the works of Donaldson \cite{D82,
D87} and Uhlenbeck-Yau \cite{UY}. It is expected that this should
also be the case for general canonical metrics in K\"ahler geometry.
For K\"ahler-Einstein metrics on Fano manifolds, obstructing sheaves
have been constructed by Nadel \cite{N} as multiplier ideal sheaves,
following ideas of Kohn \cite{JJK79} and Siu \cite{S88}.
This formulation in terms of multiplier
ideal sheaves opens up many possibilities for relations with complex
and algebraic geometry \cite{S05, DK, W}.

\medskip

The obstructing multiplier ideal sheaves are not expected to be
unique. Nadel's construction is based on the method of continuity
for solving a specific Monge-Amp\`ere equation for K\"ahler-Einstein
metrics. It has always been desirable to construct also an
obstructing multiplier ideal sheaf from the K\"ahler-Ricci
flow. The purpose of this note is to show that this can be easily
done, using the recent estimates of Kolodziej \cite{K98, K03} and
Perelman \cite{P}. In effect, Kolodziej's estimates provide a
Harnack estimate for the Monge-Amp\`ere equation, which is elliptic, and
Perelman's estimate reduces the K\"ahler-Ricci flow, which is parabolic, to
the Monge-Amp\`ere equation. Similar ideas were exploited by Tian-Zhu \cite{TZ05}
in their proof of an inequality of Harnack type for the K\"ahler-Ricci flow.

\section{The multiplier ideal sheaf}
\setcounter{equation}{0}

Let $X$ be an $n$-dimensional compact K\"ahler manifold, equipped with a K\"ahler form
$\omega_0$ with
$\mu\o_0\in c_1(X)$, where $\mu$ is a constant. The K\"ahler-Ricci
flow is the flow defined by
\bea
\label{krg}
\dot g_{\bar
kj}=-(R_{\bar kj}- \mu g_{\bar kj}),
\eea
where $g_{\bar kj}=g_{\bar kj}(t)$ is a metric evolving in time $t$
with initial value $g_{\bar
kj}(0)=\hat g_{\bar kj}$, and $R_{\bar kj}=-\pl_j\pl_{\bar k}\log\,\det
g_{\bar qp}$ is its Ricci curvature.
Since the K\"ahler-Ricci flow
preserves the K\"ahler class of the metric, we may set $g_{\bar
kj}=\hat g_{\bar kj}+\pl_j\pl_{\bar k}\phi$, and the K\"ahler-Ricci
flow can be reformulated as
\bea
\label{krphi}
\dot\phi=\log{\omega_\phi^n\over\omega_0^n}+\mu\phi-\hat f, \qquad
\phi(0)=c_0,
\eea
where we have set $\omega_\phi={i\over 2}g_{\bar
kj}dz^j\wedge d\bar z^k$, and $\hat f$ is the Ricci potential for
the metric $\hat g_{\bar kj}$, that is, the $C^\infty$ function
defined by the equation $\hat R_{\bar kj}-\mu \hat g_{\bar
kj}=\pl_j\pl_{\bar k}\hat f$, normalized by the
condition that
\bea
\int_X e^{\hat f}\omega_0^n=\int_X\omega_0^n
\equiv V.
\eea
Here and henceforth, $\hat R_{\bar kj}$ denotes the
Ricci curvature of $\hat g_{\bar kj}$, with similar conventions for
all the other curvatures of $\hat g_{\bar kj}$. The initial
potential $c_0$ is a constant, so that the initial metric coincides
with $\hat g_{\bar kj}$.
The K\"ahler-Ricci flow exists for all time $t>0$
\cite{C}, and the main issue is its convergence.
Henceforth, we shall restrict to
the case $c_1(X)>0$ of Fano manifolds unless indicated explicitly otherwise,
and set $\mu=1$.

\begin{theorem}
Let $X$ be an $n$-dimensional compact K\"ahler manifold with
$c_1(X)>0$.

\medskip

{\rm (i)} Consider the K\"ahler-Ricci flow
(\ref{krphi}) for potentials $\phi$,
with the initial value $c_0$
specified by (\ref{CT}) below.
If there exists some $p>1$ with
\bea
\label{lpbound}
{\rm sup}_{t\geq 0}\int_X e^{-p\phi}\o_0^n <\infty,
\eea
then there exists a sequence of times $t_i\to+\infty$
with $g_{\bar kj}(t_i)$ converging in $C^\infty$ to a K\"ahler-Einstein
metric. If in addition $X$ admits no non-trivial holomorphic vector field,
then the whole flows (\ref{krg}) and (\ref{krphi}) converge in $C^\infty$.

{\rm (ii)} If $X$ does not admit a K\"ahler-Einstein metric,
then for each $p>1$, there exists a function
$\psi$ which is a $L^1$ limit point of the
K\"ahler-Ricci flow (\ref{krphi}),
with the following property.
Let the multiplier ideal sheaf ${\cal J}(p\psi)$ be the sheaf
with stalk at $z$ defined by
\bea
\label{sheaf}
{\cal J}_z(p\psi)=\{f;\quad \exists U\,\ni z,\quad f\in {\cal O}(U),
\quad \int_U |f|^2 e^{-p\psi}\omega_0^n <\infty\},
\eea
where $U\subset X$ is open, and ${\cal O}(U)$ denotes the space of holomorphic functions on $U$.
Then ${\cal J}(p\psi)$ defines a proper coherent
analytic sheaf on $X$, with acyclic cohomology, i.e.,
\bea
H^q(X,K_X^{-[p]}\otimes{\cal J}(p\psi))=0,\qquad\qquad q\geq 1.
\eea
If $X$ admits a compact group $G$ of holomorphic
automorphisms, and $\hat g_{\bar kj}$ is $G$-invariant, then ${\cal J}(p\psi)$ and
the corresponding subscheme are also $G$-invariant.

\end{theorem}

\bigskip

In Part (i), once the convergence of a subsequence $g_{\bar kj}(t_i)$
has been established and $X$ is known to admit a K\"ahler-Einstein metric,
it follows from an unpublished result of Perelman that the full K\"ahler-Ricci flow
must then converge. An extension of Perelman's result to K\"ahler-Ricci
solitons is given in Tian-Zhu \cite{TZ05}.
For the sake of completeness, we have provided a short self-contained
proof of the K\"ahler-Einstein case, in our context and under the
simplifying assumption of no non-trivial holomorphic
vector fields.

\medskip

Part (ii) is of course exactly the same as in the method
of continuity for the Monge-Amp\`ere equation used by Nadel \cite{N},
in the formulation of Demailly-Koll\'ar \cite{DK}.
We divide the proof of Theorem 1 into several lemmas.

\medskip
First, we need to recall the fundamental recent result of Perelman.
Let the K\"ahler-Ricci flow be defined by (\ref{krg}), and for each
time $t$, define the Ricci potential $f$ by
\bea
\label{potential}
R_{\bar kj}-g_{\bar
kj}=\pl_j\pl_{\bar k}f, \qquad {1\over V}\int_X e^f\o_\phi^n=1.
\eea
(In particular, at time $t=0$, the Ricci potential coincides with
the function $\hat f$ defined earlier). Then Perelman \cite{P}
(see also \cite{ST}) has shown
that
\bea
\label{perelman0} {\rm sup}_{t\geq 0}(||f||_{C^0}+||\nabla
f||_{C^0}+||\Delta f||_{C^0})<+\infty,
\eea
with Laplacians and
norms taken with respect to the metric $g_{\bar kj}$.

\medskip

Next, we specify the value of the initial potential $c_0$ in
(\ref{krphi}), following Chen and Tian \cite{CT} (see
also \cite{L}). The underlying observation is that
$||\nabla\dot\phi||_{L^2}^2$, and hence the integral
\bea
\label{initial}
\int_0^\infty e^{-t}||\nabla\dot\phi||_{L^2}^2dt
\eea
does not depend on the choice of initial value $c_0$ for the
flow (\ref{krphi}). Indeed, given a flow $\phi$ with initial value
$c_0$, the function $\tilde\phi=\phi+(\tilde c_0-c_0)e^t$ satisfies
the same flow with initial value $\tilde c_0$, and hence, by
uniqueness, must coincide with the flow with initial value $\tilde
c_0$. Clearly, $\nabla\phi=\nabla\tilde\phi$, hence the assertion.
Note also that the integral in (\ref{initial}) is always finite,
in view of Perelman's estimate (\ref{perelman0}).
Following \cite{CT}, we use this common value to choose the initial value
$c_0$ in (\ref{krphi}),
\bea
\label{CT}
c_0=\int_0^\infty e^{-t}||\nabla\dot\phi||_{L^2}^2dt+{1\over
V}\int_X\hat f\o_0^n.
\eea
A specific choice of initial data is clearly necessary to discuss
the convergence of the K\"ahler-Ricci flow (\ref{krphi}) for potentials,
in view of the fact that different initial data for $\phi$
lead to flows differing by terms blowing up in time.
We will see below that the choice (\ref{CT}) is the
right choice.

\medskip

The first indication is that, with the choice
(\ref{CT}) for the initial data (\ref{krphi}), Perelman's estimate
for $h$ is equivalent to
\bea
\label{perelman}
{\rm sup}_{t\geq
0}||\dot\phi||_{C^0}<\infty.
\eea
To see this, we note that
$f+\dot\phi$ is a constant, since $\ddb(f+\dot\phi)=0$. It
suffices to show then that the average $\alpha(t)\equiv{1\over V}\int_X\dot\phi\o_\phi^n$
of $\dot\phi$ is uniformly bounded in absolute value, since $|f|$
already is, by Perelman's estimate (\ref{perelman0}). Now
differentiating the equation (\ref{krphi}) gives
$\pl_t\dot\phi=\Delta\dot\phi+\dot\phi$, and hence
\bea
\pl_t({1\over V}\int_X\dot\phi\o_\phi^n) = {1\over
V}\int_X(\Delta\dot\phi+\dot\phi)\o_\phi^n + {1\over
V}\int_X\dot\phi\Delta\dot\phi\o_\phi^n = {1\over
V}\int_X\dot\phi\o_\phi^n-||\nabla\dot\phi||_{L^2}^2.
\eea
This is a differential equation for $\alpha(t)$
which can be integrated, giving
\bea
e^{-t}\al(t)=\al(0)-\int_0^te^{-s}||\nabla\dot\phi(s)||_{L^2}^2 ds =
\int_t^\infty e^{-s}||\nabla\dot\phi(s)||_{L^2}^2ds,
\eea
in view of (\ref{CT})
and the fact that at time $t=0$, we have $\dot\phi=c_0-\hat f$.
It follows that
\bea
\label{CT1}
0\leq \al(t) = \int_t^\infty
e^{-(s-t)}||\nabla\dot\phi(s)||_{L^2}^2ds
\leq C\,\int_t^\infty
e^{-(s-t)} ds\leq C,
\eea
where we have applied Perelman's uniform
bound for $||\nabla\dot\phi||_{C^0}$.
This proves (\ref{perelman}).

\medskip
More systematically, uniform bounds for $\phi$ and $g_{\bar kj}$ are now equivalent:

\begin{lemma}
Let $X$ be a compact K\"ahler manifold, with K\"ahler form $\o_0\in c_1(X)$,
and consider the K\"ahler-Ricci flows (\ref{krg}) and (\ref{krphi}) for
$g_{\bar kj}$ and $\phi$ respectively.
Let the initial value $c_0$ for $\phi$ be
given by (\ref{CT}).
Then $||\phi||_{C^m}$ is uniformly
bounded for all $m$ if and only if $||g_{\bar k j}||_{C^m}$
is bounded for all $m$ (here norms are taken with respect to a fixed
reference metric, say $\hat g_{\bar kj}$).
The flow for $g_{\bar kj}$ converges in $C^\infty$
if and only if the flow for $\phi$ converges in $C^\infty$.

\end{lemma}

\noindent
{\it Proof of Lemma 1:} Clearly, the convergence/boundedness of the potentials $\phi$'s
implies the convergence/boundedness of the metrics $g_{\bar kj}$. Conversely,
the convergence/boundedness of the metrics implies the convergence/boundedness of $\ddb\phi$,
so it suffices to establish the convergence/boundedness
of the averages of $\phi$
with respect to the volume forms $\o_\phi^n$. Since the flow implies
\bea
{1\over V}\int_X\phi\o_\phi^n
={1\over V}\int_X \dot\phi\o_\phi^n
-
{1\over V}\int_X \log{\o_\phi^n\over\o_0^n}\o_\phi^n
+{1\over V}\int_X\hat f\o_\phi^n,
\eea
and $|\al(t)|$ is bounded in view of (\ref{CT1}), it follows that
$|{1\over V}\int_X\phi\o_\phi^n|$ is bounded
in either case. Assume now that $g_{\bar kj}$ converges. We wish to show
the convergence of ${1\over V}\int_X\phi\o_\phi^n$, and thus of $\al(t)$.

\medskip

The convergence of $g_{\bar kj}$ implies that $X$ admits a K\"ahler-Einstein metric.
By a theorem of Bando-Mabuchi \cite{BM}, the Mabuchi $K$-energy functional must
be then bounded from below.
This is well-known to imply in turn that
$||\nabla\dot\phi||_{L^2}\to 0$ as $t\to+\infty$
(see e.g. \cite{PS04} eq. (2.10) and subsequent paragraph).
But with the choice
(\ref{CT}) for initial data for (\ref{krphi}), we have the estimate
(\ref{CT1}), which implies now that $\al(t)\to 0$. Q.E.D.

\begin{lemma}
Let $X$ be a compact K\"ahler manifold, and consider the
K\"ahler-Ricci flow as defined by (\ref{krg}) and (\ref{krphi}) with
$\o_0\in c_1(X)$, and the initial value $c_0$ for $\phi$
specified by (\ref{CT}). Then for any $p>1$, we have
\bea
\label{harnack} {\rm sup}_{t\geq 0}\int_Xe^{-p\phi}\o_0^n <\infty \
\Leftrightarrow \ {\rm sup}_{t\geq 0}||\phi||_{C^0}<\infty.
\eea

\end{lemma}

\noindent
{\it Proof of Lemma 2}. This lemma is a direct consequence of the above results of Perelman
combined with results of Kolodziej. Clearly, the uniform boundedness of the $C^0$ norm of $\phi$
implies the uniform boundedness of $||e^{-\phi}||_{L^p(X)}$.
To show the converse, we consider the following Monge-Amp\`ere equation
\bea
\label{MA}
{\rm det}\,(\hat g_{\bar kj}+\pl_j\pl_{\bar k}\phi)
=
\Phi\,{\rm det}\,\hat g_{\bar kj}.
\eea
where $\Phi$ is a smooth strictly positive function. Then Kolodziej \cite{K98, K03} has shown
that, for any $p>1$, the solution $\phi$ must satisfy the a priori bound
\bea
{\rm osc}_X\phi\equiv {\rm sup}_X\phi-{\rm inf}_X\phi\ \leq\ C_p,
\eea
for some constant $C_p$ which is bounded if $||\Phi||_{L^p(X)}$ is bounded.
Now the K\"ahler-Ricci flow (\ref{krphi}) can be rewritten in the form
(\ref{MA}) with $\Phi={\rm exp}\,(\hat f-\phi+\dot\phi)$.
By Perelman's estimate (\ref{perelman}),
$||\Phi||_{L^p(X)}$ is uniformly bounded if and only if $||e^{-\phi}||_{L^p(X)}$
is uniformly bounded. Combined with Kolodziej's result, we see that the uniform
boundedness of $||e^{-\phi}||_{L^p(X)}$ implies the uniform boundedness of
${\rm osc}_X\phi$.

\medskip
To obtain a bound for $||\phi||_{C^0}$ from  ${\rm osc}\,\phi$, it suffices to
produce a lower bound for ${\rm sup}_X\phi$ and an upper bound for
${\rm inf}_X\phi$. Now, from
Perelman's estimate, we have
\bea
C_1\,e^{\hat f-\phi+\dot\phi}\o_0^n\ \leq \ e^{-\phi}\o_0^n
\leq\ C_2\, e^{\hat f-\phi+\dot\phi}\o_0^n,
\eea
and hence, integrating and recalling that
$e^{\hat f-\phi+\dot\phi}\o_0^n=\o_\phi^n$ has the same volume as $\o_0^n$,
\bea
\label{volume}
C_1\ \leq {1\over V}\int_X e^{-\phi}\o_0^n\leq\ C_2.
\eea
This implies at once that
\bea
\label{crude}
{\rm sup}_X\phi\geq -\log\,C_2,\qquad {\rm inf}_X\phi\leq -\log\,C_1.
\eea
The proof of Lemma 2 is complete. Q.E.D.

\medskip

\begin{lemma}
Let $X$ be a compact K\"ahler manifold with K\"ahler form
$\o_0$ satifying $\mu\o_0\in c_1(X)$, where $\mu$
is any constant. Let
the K\"ahler-Ricci flow be defined by (\ref{krphi}).
Then we have the a priori estimates
\bea
\label{c0}
{\rm sup}_{t\geq 0}||\phi||_{C^0}\leq A_0<\infty
\Leftrightarrow
{\rm sup}_{t\geq 0}||\phi||_{C^k}\leq A_k<\infty,
\qquad \forall\ k\in {\bf N}.
\eea

\end{lemma}

\bigskip
\noindent {\it Proof of Lemma 3} This is the parabolic analogue of
Yau's and Aubin's well-known result \cite{Y,A}, namely, that
the same statement holds for the solution $\phi$ of the elliptic
Monge-Amp\`ere equation (\ref{MA}), with the corresponding constants
$A_k$ depending on the $C^\infty$ norms of the right hand side
$\Phi$. Now the K\"ahler-Ricci flow can be rewritten in the form
(\ref{MA}), with $\Phi={\rm exp}(\hat f-\phi+\dot\phi)$. The hypothesis
$||\phi||_{C^0}\leq A_0$ implies control of $||\Phi||_{C^0}$, in
view of Perelman's estimate. However, we do not have control of all
the $C^\infty$ norms of $\Phi$, and hence Yau's a priori estimates
cannot be quoted directly.

\medskip

Thus we have to go through a full parabolic analogue of Yau's arguments,
and make sure that it goes through without any estimate on $\dot\phi$ which is not
provided by Perelman's result. The arguments here are completely parallel to
Yau's, but
we take this opportunity to present a more streamlined version.
The parabolic analogues of several key identities
are also made more explicit. They turn out to be quite simple, and may be more flexible
for future work.

\medskip
Let $\nabla$, $\Delta=\nabla^{\bar p}\nabla_{\bar p}$, $R_{\bar qp}{}^l{}_m$, etc.
and $\hat\nabla$, $\hat\Delta$, $\hat R_{\bar qp}{}^l{}_m$, etc. be the connections,
laplacians, and curvatures with respect to the metrics $g_{\bar kj}$ and $\hat g_{\bar kj}$
respectively. It is most convenient to formulate all the identities we need in terms
of the endomorphism $h=h^\al{}_\b$ defined by
\bea
h^\al{}_\b=\hat g^{\al\bar\lambda}g_{\bar\lambda\b}
\eea
For example, the difference between the connections and curvatures with respect to
$g_{\bar kj}$ and $\hat g_{\bar kj}$ can be expressed as
\bea
\label{difference}
\nabla_m V_l-\hat\nabla_m V_l&=&-V_\al(\nabla_m h\,h^{-1})^\al{}_l,
\qquad
\nabla_m V^l-\hat\nabla_m V^l=(\nabla_m h\,h^{-1})^l{}_\al V^\al
\nonumber\\
\hat R_{\bar kj}{}^\al{}_\b
-
R_{\bar kj}{}^\al{}_\b
&=&
\pl_{\bar k}(\nabla_j h\,h^{-1})^\al{}_\b
\eea
In particular, taking $V_l\to \pl_{\bar k}\pl_l\phi$, we find
\bea
\label{nabla3}
\phi_{j\bar k m}\equiv \hat\nabla_m \pl_{\bar k}\pl_j\phi
=
-g_{\bar k\al}(\nabla_m h\,h^{-1})^\al{}_j.
\eea
Henceforth, all indices are raised and lowered with respect to the metric $g_{\bar kj}$,
unless indicated explicitly otherwise. We also set
\bea
G=\log\,{\o_\phi^n\over\o_0^n}.
\eea

\medskip
\noindent
{\it Proof of the $C^2$ estimates:} The basic identity for this step is the following,
\bea
\label{identity1}
(\Delta-\pl_t)\log\,{\rm Tr}\,h
&=&
{1\over {\rm Tr}\,h}\{\hat\Delta(G-\dot\phi)-\hat R\}- {1\over {\rm Tr}\,h}
g^{p\bar q}g_{\bar mj}\hat g^{r\bar m }\hat R_{\bar qp}{}^j{}_r
\nonumber\\
&&
+
\{{\hat g^{\delta\bar k}\phi_{\g\bar k p}\phi^\g{}_\delta{}^p
\over
{\rm Tr}\,h}
-
{ g^{\delta\bar k}\pl_{\bar k}{\rm Tr}\,h\,\pl_\delta{\rm Tr}\,h
\over ({\rm Tr}\,h)^2} \}
\eea
This identity follows from another well-known identity \cite{Y}, which will also be of
later use,
\bea
\label{later}
\Delta {\rm Tr}\,h
=
\hat\Delta G -\hat R+\hat g^{\delta\bar k}\phi_{\g\bar k p}\phi^\g{}_\delta{}^p
-
g^{p\bar q}g_{\bar mj}\hat g^{r\bar m }\hat R_{\bar qp}{}^j{}_r,
\eea
and can be seen as follows:
$\Delta {\rm Tr}\,h=\bar\Delta {\rm Tr}\,h=g^{p\bar q}\nabla_{\bar q}{\rm Tr}\,
\{(\nabla_p h\,h^{-1})h\}$,
and thus
\bea
\Delta {\rm Tr}\,h
&=&
g^{p\bar q}\nabla_{\bar q}{\rm Tr}(\nabla_ph\,h^{-1})
+
g^{p\bar q}{\rm Tr}\{(\nabla_{p}h\,h^{-1})
\nabla_{\bar q}h\}
\eea
The second term on the right hand side can be recognized as
$\hat g^{\delta\bar k}\phi_{\g\bar k p}\phi^\g{}_\delta{}^p$ using (\ref{nabla3}),
while, using (\ref{difference}), the first term can be rewritten as
\bea
g^{p\bar q}\nabla_{\bar q}{\rm Tr}(\nabla_ph\,h^{-1})
=
g^{p\bar q}\hat R_{\bar qp}{}^\al{}_\b h^{\b}{}_\al-R^\al{}_\b h^\b{}_\al
=
g^{p\bar q}g_{\bar\lambda\al}\hat R_{\bar qp}{}^\al{}_\b\hat g^{\b\bar\lambda}
-R^\al{}_\b h^\b{}_\al.
\eea
But the Ricci curvature $R_{\bar\g\b}$ and be expressed in terms of $G$,
$R_{\bar\g\b}=\hat R_{\bar\g\b}-\pl_\b\pl_{\bar\g}G$. Substituting in gives
(\ref{later}). Taking the log and subtracting the simple identity
$\pl_t\,\log\,{\rm Tr}\,h=(\hat\Delta\dot\phi) ({\rm Tr}\,h)^{-1}$ gives (\ref{identity1}).

\bigskip
So far the discussion has been general. Let now $\phi$ evolve by the
K\"ahler-Ricci flow,
\bea
\dot\phi-G=\mu\phi-\hat f,
\eea
so that the term $\hat\Delta(\dot\phi-G)$ in (\ref{identity1})
can be replaced by the more tractable term $\mu\hat\Delta\phi-\hat\Delta\hat f$.
In \cite{Y}, it was shown that the expression
in brackets in (\ref{identity1}) was always non-negative, while the curvature tensor
term was bounded by
\bea
-g^{p\bar q}g_{m\bar j}\hat R^j{}_{p\bar q\al}\hat g^{\al\bar m}
=
-\sum_{i,j=1}^n {1+\phi_{\bar ii}\over 1+\phi_{\bar jj}}\hat R_{\bar ii\bar jj}
\geq -C\,({\rm Tr}\,h)\sum_j{1\over 1+\phi_{\bar jj}},
\eea
in a system of local holomorphic coordinates where both $g_{\bar kj}$ and
$\hat g_{\bar kj}$ were diagonal, and $\hat g_{\bar kj}$
was the identity matrix at a given point. Thus we have
\bea
(\Delta-\pl_t)\,\log\,{\rm Tr}\,h
\geq -\mu-C_1{1\over {\rm Tr}\,h}-C_2\sum_{j=1}^n{1\over 1+\phi_{\bar jj}}
\geq -\mu-C_3\sum_{j=1}^n{1\over 1+\phi_{\bar jj}}.
\eea
Let $A$ be any constant. Since
\bea
\Delta\phi=\sum_{j=1}^n{\phi_{\bar jj}\over 1+\phi_{\bar jj}}=n-\sum_{j=1}^n{1\over 1+\phi_{\bar jj}},
\eea
we can write
\bea
(\Delta-\pl_t)(\log{\rm Tr}\,h-A\phi\big)
\geq
C_4\dot\phi -C_5+C_6\sum_{j=1}^n{1\over 1+\phi_{\bar jj}},
\nonumber
\eea
with $A=C_4$, $C_5=\mu+An$,
and $C_6=A-C_3>0$ for $A$ large enough. In view of Perelman's estimate (\ref{perelman}), we conclude
\bea
(\Delta-\pl_t)(\log\,{\rm Tr}\,h-A\phi)
\geq
-C_7+C_6\sum_{j=1}^n{1\over 1+\phi_{\bar jj}}.
\eea
Let now $[0,T]$ be any time interval, and $(z_0,t_0)$ a point in $X\times [0,T]$ where
the function $\log\,{\rm Tr}\,h-A\phi$ attains its maximum. If this point
is not at time $t=0$, then the left hand side
of the above equation is $\leq 0$, and we obtain the estimate
\bea
{1\over 1+\phi_{\bar jj}}\leq C_8,
\qquad 1\leq j\leq n.
\eea
But then, at the point $(z_0,t_0)$,
\bea
{\rm Tr}\,h
&=&
{\rm Tr}\,h({{\rm det}\, \hat g_{\bar kj}\over{\rm det}\,g_{\bar kj}})e^{\hat f-\phi+\dot\phi}
=
e^{\hat f-\phi+\dot\phi}\sum_{i=1}^n(1+\phi_{\bar ii})\prod_{j=1}^n{1\over 1+\phi_{\bar jj}}
\nonumber\\
&=&e^{\hat f-\phi+\dot\phi}\sum_{i=1}^n\prod_{j\not=i}{1\over 1+\phi_{\bar jj}}
\leq C_9,
\eea
using the boundedness of $||\phi||_{C^0}$ and again Perelman's estimate. But now
we have
\bea
{\rm sup}_{X\times[0,T]} {\rm Tr}\,h\leq e^{A||\phi||_{C^0}}{\rm exp}{\big(\log\,
{\rm Tr}\,h-A\phi\big)(z_0,t_0)}
\leq C_{10}.
\eea
Since $T$ is arbitrary, this establishes the boundedness of the trace of
$\hat g_{\bar kj}+\pl_j\pl_{\bar k}\phi$, and since the matrix is positive, of all
its entries. The proof of the $C^2$ estimate is complete.

\medskip
\noindent {\it Proof of the $C^3$ estimates:} this step was
established in \cite{C} when $c_1(X)=0$ or $c_1(X)<0$, and in
\cite{P} for $c_1(X)>0$. We shall give below a simpler proof for all
cases with completely explicit formulas.

\medskip

The main ingredient is a parabolic analogue of the Yau, Aubin, and Calabi identities
for the third derivatives of the Monge-Amp\`ere equation. In their case,
the Ricci curvature is pre-assigned and hence all its derivatives can be
controlled. In the present case, we cannot control as yet the derivatives
of the Ricci curvature, and it is crucial that they cancel out in the desired
identity. We show this by a completely explicit formula,
the main technical innovation being the use of the endomorphism
$h^\al{}_\b$ instead of the potential $\phi$ itself.
The squared terms in the Calabi identity arise naturally as
the familiar squared terms in a formula of Bochner-Kodaira type.

\medskip
Let $S$ be defined as in \cite{Y} by
\bea
S=
g^{j\bar r} g^{s\bar k}g^{m\bar t}  \phi_{j\bar k m}\phi_{\bar r s\bar t}
\eea
In terms of $h^\al{}_\b$, $S$ is just the square of the $L^2$ norm of the $g_{\bar kj}$ connection,
\bea
S=
g^{m\bar\g}g_{\bar\mu \b}g^{l\bar\al}
(\nabla_mh\,h^{-1})^\b{}_l\overline{(\nabla_\g h\,h^{-1})^\mu{}_\al}
=
|\nabla h\,h^{-1}|^2
\eea
and its Laplacian leads immediately to a formula of Bochner-Kodaira type,
\bea
\Delta S&=&
g^{m\bar\g}g_{\bar\mu \b}g^{l\bar\al}
(\ \Delta(\nabla_mh\,h^{-1})^\b{}_l\overline{(\nabla_\g h\,h^{-1})^\mu{}_\al}
+(\nabla_mh\,h^{-1})^\b{}_l\overline{\bar\Delta(\nabla_\g h\,h^{-1})^\mu{}_\al}\ )
\nonumber\\
&&
+
|\bar\nabla(\nabla h\,h^{-1})|^2+|\nabla(\nabla h\,h^{-1})|^2
\eea
where, more explicitly,
$|\bar\nabla(\nabla h\,h^{-1})|^2=
g^{q\bar p}g^{j\bar m}
g^{\b\bar\delta}g_{\bar\g \al}\nabla_{\bar p}(\nabla_jh\,h^{-1})^\al{}_\b
\overline{\nabla_{\bar q}(\nabla_mh\,h^{-1})^\g{}_\delta}$, etc.
The relation between $\bar\Delta$ and $\Delta$
follows from commuting the $\nabla_q$ and the $\nabla_{\bar p}$ derivatives,
\bea
(\bar\Delta(\nabla_j h\,h^{-1}))^\g{}_\al
&=&
(\Delta(\nabla_j h\,h^{-1}))^\g{}_\al
-
R^\g{}_\mu (\nabla_\g h\,h^{-1})^\mu{}_\al
+
R^\mu{}_\al(\nabla_j h\,h^{-1})^\g{}_\mu
\nonumber\\
&&
+
R^\mu{}_j(\nabla_\mu h\,h^{-1})^\g{}_\al
\eea
Thus we have
\bea
\label{explicit0}
\Delta S&=&
g^{m\bar\g}g_{\bar\mu \b}g^{l\bar\al}
(\ \Delta(\nabla_mh\,h^{-1})^\b{}_l\}\overline{(\nabla_\g h\,h^{-1})^\mu{}_\al}
+
(\nabla_mh\,h^{-1})^\b{}_l\overline{\Delta(\nabla_\g h\,h^{-1})^\mu{}_\al}\ )
\nonumber\\
&&+
|\bar\nabla(\nabla h\,h^{-1})|^2+|\nabla(\nabla h\,h^{-1})|^2
\nonumber\\
&&
+(\nabla_m h\,h^{-1})^\b{}_l
(\ g^{m\bar\g}g_{\bar\mu\b}R^{l\bar\rho}\overline{(\nabla_\g h\,h^{-1})^\mu{}_\rho}
-
g^{m\bar \g}R_{\bar\rho\b}g^{l\bar\al}\overline{(\nabla_\g h\,h^{-1})^\rho{}_\al}
\nonumber\\
&&
\qquad\qquad\qquad \qquad
+
R^{m\bar\rho}g_{\bar\mu\b}g^{l\bar\al}\overline{(\nabla_\rho h\,h^{-1})^\mu{}_\al}
\ )
\eea
In the case of the elliptic Monge-Amp\`ere equation, this equation suffices already to
establish the desired inequality $\Delta S\geq -C_1 S-C_2$. This is because the
Ricci tensor $R_{\bar\al\b}$ is known in that case, and
the Laplacian of $\nabla h\,h^{-1}$ can be readily reduced to $\nabla R_{\bar\al\b}$,
using (\ref{difference}) and the Bianchi identity,
\bea
\Delta (\nabla_j h\,h^{-1})^l{}_m
=\nabla^{\bar p}\pl_{\bar p}(\nabla_jh\,h^{-1})
=
-\nabla^{\bar p}R_{\bar pj}{}^l{}_m
-
\nabla^{\bar p}\hat R_{\bar kpj}{}^l{}_m
=
-\nabla_jR^l{}_m
+
\nabla^{\bar p}\hat R_{\bar p j}{}^l{}_m
\nonumber
\eea
Since the connection $\nabla^{\bar p}$ is manifestly $O(\nabla h\,h^{-1})=O(\sqrt S)$,
the desired lower bound follows at once. The full expression for $\Delta S$
in terms of $R_{\bar\al\b}$ may also be of interest,
\bea
\label{explicit1}
\Delta S
&=&
-g^{m\bar\g}g_{\bar\mu\b}g^{l\bar\al} (\ \nabla_mR^\b{}_l\overline{(\nabla_\g h\,h^{-1})^\mu{}_\al}
+
(\nabla_m h\,h^{-1})^\b{}_l\overline{\nabla_\g R^\mu{}_\al}\ )
\nonumber\\
&&
+|\bar\nabla(\nabla h\,h^{-1})|^2+|\nabla(\nabla h\,h^{-1})|^2
\nonumber\\
&&
+(\nabla_m h\,h^{-1})^\b{}_l
(\ g^{m\bar\g}g_{\bar\mu\b}R^{l\bar\rho}\overline{(\nabla_\g h\,h^{-1})^\mu{}_\rho}
-
g^{m\bar g}R_{\bar\rho\b}g^{l\bar\al}\overline{(\nabla_\g h\,h^{-1})^\rho{}_\al}
\nonumber\\
&&
\qquad\qquad\qquad \qquad
+
R^{m\bar\rho}g_{\bar\mu\b}g^{l\bar\al}\overline{(\nabla_\rho h\,h^{-1})^\mu{}_\al}
\ )
\nonumber\\
&&
+
g^{m\bar\g}g_{\bar\mu\b}g^{l\bar\al}
(\ \nabla^{\bar p}\hat R_{\bar p m}{}^\b{}_l\overline{(\nabla_\g h\,h^{-1})^\mu{}_\al}
+
(\nabla_m h\,h^{-1})^\b{}_l\overline{\nabla^{\bar p}\hat R_{\bar p \g}{}^\mu{}_\al}\ ).
\eea

\bigskip

In the parabolic case, we do not have control of
$R_{\bar\al\b}$ and its derivatives,
and need to eliminate these terms using the time derivative $\dot S$
of $S$. We begin by giving a general formula for $\dot S$ in terms of $h^{-1}\dot h$,
so that it is valid for all evolutions.
First, note that the derivatives of the connection are given by
\bea
(\nabla_j h)^{\dot{}}
=\nabla_j(h^{-1}\dot h)h+(\nabla_jh)(h^{-1}\dot h),
\qquad
(\nabla_jh\,h^{-1})^{\dot{}}
=\nabla_j(h^{-1}\dot h).
\eea
Indeed, note that
$\dot g_{\bar\al\b}=g_{0\bar\al\mu}\dot h^{\mu}{}_\b=g_{\bar\al\nu}(h^{-1}\dot h)^\nu{}_\b
=
(gh^{-1}\dot h)_{\bar\al\b}$, and
$(g^{\b\bar b})^{\dot{}}=-(\dot g)^{\b\bar b}
=
-(h^{-1}\dot h)^\b{}_\nu g^{\nu\bar b}
=
-(h^{-1}\dot h g^{-1})^{\b\bar b}$.
Writing
\bea
\{(\nabla_jh)^{\dot{}}\}^p{}_q
=
g_{\bar b q}g^{p\bar a}\pl_j(g^{\b\bar b}g_{\bar a\al}h^\al{}_\b)
=
\{g^{-1}\pl_j(ghg^{-1})g\}^p{}_q
\eea
differentiating with respect to time, and substituting in the
preceding formulas for $\dot g_{\bar \al\b}$ and $(\bar g^{\bar\al\b})^{\dot{}}$ gives
at once
$(\nabla_j h)^{\dot{}}
=-h^{-1}\dot h\nabla_j h+(\nabla_jh)(h^{-1}\dot h)+\nabla_j(h^{-1}\dot h h)$,
from which the desired formula for $(\nabla_j h)^{\dot{}}$ follows.
The formula for $(\nabla_jh h^{-1})^{\dot{}}$ is a simple consequence of
the one for $(\nabla_j h)^{\dot{}}$. Next, differentiating $S$ gives
\bea
\dot S&=&
+g^{m\bar\g}g_{\bar\mu\b}g^{l\bar\al}
(\ \pl_t
((\nabla_mh\,h^{-1})^\b{}_l)
\overline{(\nabla_\g h\,h^{-1})^\mu{}_\al}+
(\nabla_mh\,h^{-1})^\b{}_l)
\overline{\pl_t(\nabla_\g h\,h^{-1})^\mu{}_\al}
\ )
\nonumber\\
&&
- (\nabla_mh\,h^{-1})^\b{}_l (\
(h^{-1}\dot h)^{m\bar\g}g_{\bar\mu\b}g^{l\bar\al}
\overline{(\nabla_\g h\,h^{-1})^\mu{}_\al}
+
g^{m\bar\g}(h^{-1}\dot h)_{\bar \mu\b}g^{l\bar\al}
\overline{(\nabla_\g h\,h^{-1})^\mu{}_\al}
\nonumber\\
&&
\qquad\qquad\qquad\qquad
- g^{m\bar\g}g_{\bar\mu\b}(h^{-1}\dot h)^{l\bar\al}
\overline{(\nabla_\g
h\,h^{-1})^\mu{}_\al}\ )
\eea
Combining with (\ref{explicit0}), we obtain the following general heat equation,
\bea
(\Delta-\pl_t)S
&=& |\bar\nabla(\nabla h\,h^{-1})|^2+|\nabla(\nabla h\,h^{-1})|^2
\nonumber\\
&&
+
g^{m\bar\g}g_{\bar\mu\b}g^{l\bar\al}
\big\{
(\Delta-\pl_t)(\nabla_mh\,h^{-1})^\b{}_l
\overline{(\nabla_\g h\,h^{-1})^\mu{}_\al}
\nonumber\\
&&
+g^{m\bar\g}g_{\bar\mu\b}g^{l\bar\al}
\big\{
(\nabla_mh\,h^{-1})^\b{}_l
\overline{(\Delta-\pl_t)(\nabla_\g h\,h^{-1})^\mu{}_\al}
\nonumber\\
&&
+
\bigg\{(h^{-1}\dot h+R)^{m\bar\g}g_{\bar\mu\b}g^{l\bar\g}
-
g^{m\bar\g}(h^{-1}\dot h+R)_{\bar\mu\b}g^{l\bar\al}
+
g^{m\bar\g}g_{\bar\mu\b}(h^{-1}\dot h+R)^{l\bar\al}\bigg\}
\nonumber\\
&&
\qquad\qquad\times
(\nabla_mh\,h^{-1})^\b{}_l
\overline{(\nabla_\g h\,h^{-1})^\mu{}_\al}
\eea
We can now specialize to the K\"ahler-Ricci flow, where
\bea
(h^{-1}\dot h)^\b{}_l=- (R^\b{}_l-\mu\,\delta^\b{}_l),
\eea
and hence, using again (\ref{difference}),
\bea
(h^{-1}\dot h+R)^\b{}_\l
&=&\mu\,\delta^\b{}_\l
\nonumber\\
(\Delta-\pl_t)(\nabla_jh\,h^{-1})^l{}_m
&=&
\nabla^{\bar q}\hat R_{\bar q j}{}^l{}_m.
\eea
Substituting this in the previous formula for $(\Delta-\pl_t)S$, we obtain
the following simple and completely explicit parabolic analogue
of the $C^3$ identity of Yau, Aubin, and Calabi,
\bea
(\Delta-\pl_t)S
&=&
|\bar\nabla(\nabla h\,h^{-1})|^2+|\nabla(\nabla h\,h^{-1})|^2
+\mu|\nabla h\,h^{-1}|^2
\nonumber\\
&&
+
g^{m\bar\g}\nabla^{\bar q}\hat R_{\bar q m}{}^\b{}_l
\overline{(\nabla_\g h\,h^{-1})_{\bar\b}{}^{\bar l}}
+
g^{m\bar \g}(\nabla_m h\,h^{-1})_{\bar\mu}{}^{\bar\al}
\overline{\nabla^{\bar q}\hat R_{\bar q\g}{}^\mu{}_\al}
\eea

Note that the terms in $R_{\bar kj}$ and its derivatives have cancelled out.
Since $\hat R_{\bar q m}{}^l{}_\b$ is a fixed tensor,
we obtain immediately the estimate
\bea
(\Delta-\pl_t)S
\geq
|\bar\nabla(\nabla h\,h^{-1})|^2+|\nabla(\nabla h\,h^{-1})|^2
-C_1\,S-C_2.
\eea
The $C^3$ estimates for $\phi$ can now be established by the standard arguments:
by the $C^2$ estimates, the metric $g_{\bar kj}$ is known to be equivalent
to $\hat g_{\bar kj}$. Thus $S$ is of the same size as
the expression $\hat g^{\delta\bar k}\phi_{\g\bar k p}\phi^\g{}_\delta{}^p$.
In view of the identity (\ref{later}), we have
for
$A$
sufficiently large,
\bea
(\Delta-\pl_t)(S+A\hat\Delta\phi)
\geq
C_3 S-C_4,
\eea
with $C_3>0$. The maximum principle implies now that $S$ is bounded by a fixed
positive constant. This completes the proof of the $C^3$ estimates.

\bigskip
The remaining part of the proof of Lemma 3 is standard: the uniform $C^0$ and $C^2$
estimates for $\phi$ imply uniform $C^1$ estimates, and together with
the uniform $C^3$ estimates, we can deduce that the flow (\ref{krphi})
is a parabolic equation with uniform $C^1$ coefficients. The general theory
of parabolic PDE's can then be applied to give uniform $C^k$ estimates for all orders $k$.
Q.E.D.

\bigskip
\noindent
{\it Proof of Theorem 1, part (i):}
By Lemmas 2 and 3,
the uniform estimate (\ref{lpbound}) implies uniform estimates for
$||\phi||_{C^m}$ for each $m\in{\bf N}$. The Arzela-Ascoli theorem implies
the existence of times $t_i\to+\infty$ with $\phi(t_i)$ and
$g_{\bar kj}(t_i)$ converging in $C^\infty$.

\medskip
To show that the limit of the subsequence $g_{\bar kj}(t_i)$
is a K\"ahler-Einstein metric,
we observe that the uniform boundedness of $||\phi||_{C^m}$
for each $m\in{\bf N}$ implies that the metrics $g_{\bar kj}$
are all equivalent along the K\"ahler-Ricci flow,
and that their curvature tensors and their derivatives are
all uniformly bounded.
This is easily seen to imply that the Mabuchi $K$-energy $\nu_{\o_0}(\phi)$
is bounded along the flow, since $\nu_{\o_0}(\phi)$
can be written explicitly as
\bea
\label{functional}
\nu_{\o_0}(\phi) =
{1\over V}
\int_X\{\log ({\o_\phi^n\over\o_0^n})\o_\phi^n
-
\phi(Ric(\o_0)\sum_{i=0}^{n-1}\o_0^i\o_\phi^{n-1-i}
-
{n\over n+1}\sum_{i=0}^n\o_0^i\o_\phi^{n-i})\}.
\eea
(see also \cite{DT,Ch} for alternative expressions).
Next, the lower bound for the Mabuchi energy and the
uniform boundedness of the curvature tensor together imply
(\cite{PS04}, Theorem 1)
\bea
||R_{\bar kj}-g_{\bar kj}||_{L^2}\to 0,
\qquad t\to+\infty,
\eea
where the $L^2$-norm is with respect to $g_{\bar kj}$.
Since $g_{\bar kj}$ is uniformly equivalent to $\hat g_{\bar kj}$,
this holds also for the $L^2$ norm with respect to $\hat g_{\bar kj}$.
Returning to the subsequence $g_{\bar kj}(t_i)$ converging in $C^\infty$,
the limit is then a smooth metric $g_{\bar kj}^\infty$
with $R_{\bar kj}^\infty-g_{\bar kj}^\infty=0$,
as was to be proved.

\medskip

To show the full convergence of the K\"ahler-Ricci flow when $X$
does not admit non-trivial holomorphic vector fields, we show first
that, under the assumption of uniform bounds for $||\phi||_{C^m}$
(and hence for $||g_{\bar kj}||_{C^m}$) for each $m\in{\bf N}$, the
functions $\dot\phi$ converge in $C^\infty$. We already know that
$\pl_j\pl_{\bar k}\dot\phi=-(R_{\bar kj}-g_{\bar kj})\to 0$
in all Sobolev norms if the curvature is
bounded and the Mabuchi energy is bounded from below (\cite{PS04},
Theorem 1).
Since the metrics $g_{\bar kj}$ are all equivalent, and
all their derivatives are bounded, it follows that
$\ddb\dot\phi\to 0$ in $C^\infty$. Next,
from the lower boundedness of the Mabuchi
energy and the proof of Lemma 1, the averages
$\al(t)={1\over V}\int_X\dot\phi\o_\phi^n$
of $\dot\phi$ converge to $0$ as $t\to\infty$.
To deduce the convergence of $\dot\phi$ from the convergence of
its averages,
we write for any constants $\delta,A$,
\bea
\dot\phi(z)={1\over V}\int_X
(G(z,w)+A)(-\Delta\dot\phi(w)+\delta)\o_\phi^n -\delta\,A+{1\over
V}\int_X\dot\phi\o_\phi^n,
\eea
where $G(z,w)$ is the Green's
function with respect to $g_{\bar kj}$. Since $R_{\bar kj}\to
g_{\bar kj}$ in $C^\infty$, by a theorem of Cheng-Li \cite{CL}, it
follows that we can choose a fixed $A>0$ so that $G(z,w)+A>0$ for
all $t$. For any $\epsilon>0$, we can choose $T$ large enough so
that $0\leq {1\over V}\int_X\dot\phi\o_\phi^n<{\epsilon\over 2}$ and
$-\Delta\dot\phi+{\epsilon\over 2A}>0$ for all $t>T$. It follows from
the preceding identity with $\delta={\epsilon\over 2A}$ that
\bea
{\rm inf}_X\dot\phi\ \geq -{\epsilon\over 2}, \qquad t>T.
\eea
The same argument
with $\dot\phi$ replaced by $-\dot\phi$ gives the bound ${\rm
sup}_X\dot\phi\leq \epsilon$, and thus we have shown that
$\dot\phi\to 0$ in $C^0$. Together with the convergence of
$\ddb\dot\phi$ to $0$ in $C^\infty$, this implies that $\dot\phi\to
0$ in $C^\infty$.

\medskip

Next, note that the operator $\Delta(t)+1$
is uniformly bounded away from $0$, where $\Delta(t)=\nabla^{\bar p}\nabla_{\bar p}$ is
the scalar
Laplacian defined by the evolving metric $g_{\bar kj}(t)$.
Indeed, assume otherwise. Then there exists a subsequence $t_i\to+\infty$
with Laplacians $\Delta(t_i)+1$ admitting eigenvalues $\lambda(t_i)\to 0$.
By going to a subsequence, we may assume that the corresponding metrics
$g_{\bar kj}(t_i)$ converge in $C^\infty$ to a metric $g_{\bar kj}^\infty$
which is K\"ahler-Einstein,
by the preceding discussion. The Laplacian $\Delta(\infty)+1$
admits then a zero eigenvalue. The corresponding eigenfunction $u(z)$
defines then a non-trivial holomorphic vector field $V^j=g^{\infty j\bar k}\pl_{\bar k}u$,
which contradicts our assumption.

\medskip
Fix now $t$, and consider the equation in $\psi$
\bea
\label{linearized}
\log{(\o_{\phi(t)}+{i\over 2}\ddb\psi)^n\over\o_{\phi(t)}^n}
+
\psi=h
\eea
The linearization of the left hand side at $\psi=0$ is $\Delta(t)+1$.
Since this operator is invertible, with uniform bounds in $t$,
and since the higher order derivatives of the left hand side,
viewed as a functional in $\psi$, are also uniformly bounded,
it follows that there exist constants $\epsilon_m>0$, $A_m<\infty$,
independent of $t$, so
that
\bea
\label{psi}
||\psi||_{C^{m+1}}\leq A_m\,||h||_{C^m},
\qquad
\qquad{\rm when}\ ||h||_{C^m}<\epsilon_m.
\eea

\medskip
Finally, let $t,s>>1$. Then the function $\psi\equiv \phi(s)-\phi(t)$
satisfies the equation (\ref{linearized}), with right hand side
$h$ given by
\bea
h=\dot \phi(t)-\dot\phi(s).
\eea
For any $\epsilon>0$, we have
$||\dot\phi(t)-\dot\phi(s)||_{C^m}<\epsilon$
for $s,t$ large enough
since $\dot\phi\to 0$ in $C^m$.
It follows from (\ref{psi})
that $||\phi(s)-\phi(t)||_{C^{m+1}}\leq A_m\epsilon$.
This establishes the convergence of the K\"ahler-Ricci flow
(\ref{krphi}), and hence also of (\ref{krg}). Q.E.D.

\bigskip
\noindent
{\it Proof of Theorem 1, part (ii):}
We follow here closely the arguments of Nadel \cite{N} and Demailly-Koll\'ar \cite{DK}.
If $X$ does not admit a K\"ahler-Einstein metric, then by Part (i) of the theorem,
for any $p>1$, there must exist a sequence $\phi(t_i)$ in the K\"ahler-Ricci flow with
\bea
{\rm lim}_{i\to\infty}\int_X e^{-p\phi(t_i)}\o_0^n=\infty.
\eea
Let $\psi$ be an $L^1$ limit point of the sequence $\phi(t_i)$.
Then by the Demailly-Koll\'ar theorem
on the semi-continuity of complex singularity exponents
(\cite{DK}, Main Theorem), $||e^{-\psi}||_{L^p(X)}=+\infty$,
and hence the multiplier ideal sheaf
${\cal J}(p\psi)$ is non-trivial.
Equivalently, the corresponding subscheme of structure
sheaf ${\cal O}_X/{\cal J}(p\psi)$ is non-empty.
Since $p\psi$ is strictly plurisubharmonic with respect to
$([p]+1)\o_0$,
by the theorem of Nadel \cite{N}, in the formulation
of Demailly and Koll\'ar (\cite{DK}, Theorem 4.1 and Corollary 6.6), the multiplier ideal sheaf
${\cal J}(p\psi)$ is a coherent analytic sheaf in ${\cal O}_X$ with
$K_X^{-[p]}\otimes {\cal J}(p\psi)$ having acyclic cohomology,
with the $G$-invariance property stated if the metric $\hat g_{\bar kj}$ is $G$-symmetric.
Q.E.D.

\section{Remarks}
\setcounter{equation}{0}

We conclude with a few simple remarks.

\medskip

$\bullet$ The estimate (\ref{volume}) implies that the
functionals $F_{\o_0}(\phi)$ and $F_{\o_0}^0(\phi)$
defined by
\bea
\label{functionals} F_{\o_0}(\phi) &=&
F_{\o_0}^0(\phi)
- \log({1\over V}\int_X e^{\hat f-\phi}\o_0^n)
\nonumber\\
F_{\o_0}^0(\phi)&=&
J_{\o_0}(\phi)-{1\over
V}\int_X\phi\o_0^n
\eea
where $J_{\o_0}(\phi)={1\over 2V}\int_X\phi(\o_0^n-\o_\phi^n)$ is the
Aubin-Yau functional,
are bounded by one another along the K\"ahler-Ricci flow, up to additive constants,
\bea
-C_3+F_{\o_0}^0(\phi)\ \leq
\ F_{\o_0}(\phi)\ \leq\ C_4+F_{\o_0}^0(\phi).
\eea
This may be of some practical use, since while the functional $F_{\o_0}(\phi)$ decreases along
the K\"ahler-Ricci flow, it is the functional $F_{\o_0}^0(\phi)$ which is more directly linked to
Chow-Mumford stability (see e.g. \cite{Z, PS02, Paul}).

\medskip
$\bullet$ Since the functional $F_{\o_0}(\phi)$ decreases along the K\"ahler-Ricci flow,
$F_{\o_0}(\phi)$ and hence $F_{\o_0}^0(\phi)$ are bounded from above. Thus, along the
K\"ahler-Ricci flow, we have
\bea
\label{inequalityJ1}
J_{\o_0}(\phi)\leq {1\over V}\int_X\phi\,\o_0^n+C.
\eea

\medskip
$\bullet$ The Harnack inequality along the K\"ahler-Ricci flow proved by Tian-Zhu \cite{TZ05}
is the following,
\bea
\label{harnackkr}
{\rm osc}_X\phi \leq C_5\, J_{\o_0}(\phi)^{n+\delta}+C_6,
\eea
where $\delta$ is any positive number less or equal to $1$. As shown in \cite{TZ05},
Perelman's bound for $||\dot\phi||_{C^0}$ allows arguments similar to the proof of
the Harnack inequality for the Monge-Amp\`ere equation $\o_\phi^n=e^{F-t\phi}\o_0^n$
\cite{S88, T87} to take over and establish (\ref{harnackkr}).

\medskip
$\bullet$ The $C^0$ estimate $||\phi||_{C^0}\leq C$ along the K\"ahler-Ricci flow is now
seen to be equivalent to the following estimate on the averages of $\phi$,
\bea
\label{average}
{1\over V}\int_X\phi\,\o_0^n\leq C.
\eea
Indeed, this inequality together with the inequalities (\ref{inequalityJ1}) and (\ref{harnackkr})
imply that ${\rm osc}_X\phi$ is uniformly bounded, which implies in turn that $||\phi||_{C^0}$
is uniformly bounded, in view of the crude bounds (\ref{crude}) for ${\rm sup}_X\phi$
and ${\rm inf}_X\phi$.
The condition $\o_\phi^n\geq C\,\o_0^n$ along the K\"ahler-Ricci flow
introduced by Pali \cite{Pa} can be interpreted in this light: by Perelman's estimate,
this condition is equivalent to $e^{-\phi}\geq C_1>0$, and hence ${\rm sup}_X\phi\leq C_2$,
which implies (\ref{average}).

\medskip
$\bullet$ The application of Kolodziej's theorem to the K\"ahler-Ricci flow confirms that
in his theorem, $p$ cannot be taken to be $1$: indeed, $e^{-\phi}$ is uniformly in $L^1(X)$
by the estimate (\ref{volume}), so that if the theorem holds for $p=1$, it would follow
from the above arguments that the K\"ahler-Ricci flow always admits
a subsequence converging to a K\"ahler-Einstein
metric on a K\"ahler manifold $X$ with $c_1(X)>0$, which is
known not to be the case.

\medskip
$\bullet$ In the applications of multiplier ideal sheaves to the existence of K\"ahler-Einstein
metrics, only the non-existence of non-trivial acyclic multiplier ideal sheaves
has played a role so far, and not the value of $p$. It is conceivable that the value of
$p$ can carry information. If so, then the multiplier ideal sheaves arising from
the K\"ahler-Ricci flow may carry more information than their counterparts from
the method of continuity: indeed, as noted above, the condition $p>1$ for
the K\"ahler-Ricci flow is sharp, while this is not known for the condition
${n\over n+1}<p<1$ for the method of continuity.

\bigskip

\bigskip
\noindent
{\it Acknowledgements.} We would like to thank Gabor
Szekelyhidi for pointing out to us an oversight in the application
of the Nadel vanishing theorem in the original version of
this paper.

\bigskip
$^*$ Department of Mathematics

Columbia University, New York, NY 10027\\

\v

$^{\dagger}$ Department of Mathematics

Rutgers University, Newark, NJ 07102

\bigskip

\end{document}